\documentclass[12pt]{article}

\usepackage{amsmath}
\usepackage{amsfonts}
\usepackage{amsthm}     
\usepackage{amssymb}  
\usepackage{latexsym}
\usepackage{euscript}    
\usepackage{txfonts}

\newcommand{\mcal}[1]{{\mathcal {#1}}} 

\theoremstyle{definition}
\theoremstyle{definition}
\theoremstyle{definition}
\theoremstyle{definition}
\theoremstyle{definition}
\theoremstyle{definition}

\theoremstyle{plain}\newtheorem*{lemma*}{Theorem}

\theoremstyle{plain}\newtheorem*{theorem*}{Theorem}
\theoremstyle{plain}\newtheorem*{corollary*}{Corollary}
\theoremstyle{remark}\newtheorem*{remark*}{Remark}
\theoremstyle{remark}\newtheorem*{remarks*}{Remarks}
\theoremstyle{definition}\newtheorem*{conjecture*}{Conjecture}
\theoremstyle{definition}\newtheorem*{definition*}{Definition}
\theoremstyle{definition}\newtheorem*{example*}{Example}

\theoremstyle{definition}\newtheorem*{question*}{Question}
\theoremstyle{definition}\newtheorem*{questions*}{Questions}

\theoremstyle{definition}\newtheorem*{hypothesis*}{Hypothesis}

\def\qed{\ifhmode\unskip\nobreak\fi\ifmmode\ifinner\else\hskip5pt\fi\fi
 \hfill\hbox{\hskip5pt\vrule width4pt height6pt depth1.5pt\hskip1pt}}
 
\newcommand{\Int}{\ensuremath{\operatorname{\mathsf{Int}}}}

\newcommand{\dist}{\ensuremath{\operatorname{\mathsf{dist}}}}

\newcommand{\dimply}{\ensuremath{\:\Longleftrightarrow\:}}

\newcommand{\R}[1]{\ensuremath{{\mathbb R}^{#1}}} 

\newcommand{\pde}[2]
	{\ensuremath{\frac{\partial #1}{\partial #2}}}	
\newcommand{\ode}[2]{\ensuremath{\frac{d#1}{d#2}}}    	

\def\d#1dt{\frac{d#1}{dt}}    

\newcommand{\del}{\ensuremath{\delta}}

\newcommand{\co}{\colon\thinspace} 

\newcommand{\sm}[1]{\ensuremath{\setminus #1}}

 \def\Rnp{{\R n}_+}

\reversemarginpar   

\begin{document}

\title{\bf On the nonchaotic nature of monotone dynamical systems}
 \author{ Morris W. Hirsch\thanks{I am grateful to Professor Robert
     Devaney for helpful information and comments.}\\ Department of
   Mathematics\\ University of Wisconsin at Madison\\ University of
   California at Berkeley}
\maketitle

\begin{abstract} 
  \noindent Two commonly used types of dynamics, chaotic and
  monotone, are compared.  It is shown that monotone maps in strongly
  ordered spaces do not have chaotic attracting sets.
  
\end{abstract}

\tableofcontents

\section{Introduction}   

This article contrasts two protean types of dynamical systems, {\em
  chaotic} and {\em monotone}.  Both types occur frequently in
mathematical models of applied fields, including Biology, Chemistry,
Economics and Physics. 

\smallskip
Topological spaces, denoted by capital Roman letters, have metrics
denoted by $d$.  The interior of $P\subset Y$ is $\,\Int \,P$.  The
distance from $y\in Y$ to $Q\subset Y$ is
\[\dist\,(y,Q):=\inf\{d(y,q)\co q\in Q\}.\]

Maps are assumed to be continuous.  The {\em orbit} of $y$ under $T\co
Y\to Y$ is the set
$\{y,Ty, T^2x, \dots\}$.
If this orbit finite,
$y$ and its orbit are {\em periodic}.

A set $A\subset Y$ {\em attracts} the point $y \in Y$ if

\begin{equation}\label{eq:1}
 \lim_{k\to\infty}d(T^ky, A) =0.
\end{equation}
We call $A$  {\em attracting} for $T$ provided:

   \begin{itemize}

   \item  $TA\subset A$,

   \item $A$ is compact and nonempty,

   \item $A$ attracts every point in some open neighborhood  $W$ of $A$
     such that $TW\subset W$. 
   \end{itemize}
We call   $A$ is an
{\em attractor} for $T$ when the limit in (\ref{eq:1}) is uniform in
$W$.
The  union of the subbasins of $A$ is  the {\em basin} of $A$,  an
open set invariant under $T$.  An attractor is {\em global} if its basin is all of $Y$.

\subsection{Chaotic dynamics.} 
The hallmark of a chaotic system is  ``sensitive dependence on
initial conditions:'' there exists  $\del >0$ such that for
any two
initial states $x(0), y(0)$, at some time $t>0$ the distance from
$x(t)$ to $y(t)$ exceeds $\del$.  If a system with this feature models
a natural system (e.g., weather, economy, ecology, gene system, a
disease) then it cannot be used to make accurate long-term
predictions.

This was discovered by the meterologist {\sc Edward Lorenz} in his
seminal 1963 article, ``Deterministic Non-periodic Flow''
\cite{Lorenz63}.  After drastically simplifying standard equations for
fluid flow, Lorenz arrived at the system of differential equations

\begin{equation*}
\begin{split}
  \dot x &= 10(y - z),\\
  \dot y &=  28x -y   - xz,\\
  \dot z &= xy - (8/3)z.
\end{split}
\end{equation*}
Despite its simple algebraic form, Lorenz found a
disturbing feature in his numerical solution:
\begin{quote} \dots two states differing by imperceptible amounts may
  eventually evolve into two considerably different states. If, then,
  there is any error whatever in observing the present state--- and in
  any real system such errors seem inevitable--- an acceptable
  prediction of an instantaneous state in the distant future may well
  be impossible.
\end{quote}
Lorenz's extensive computations convincingly illustrate this
phenomenon--- an
unexpected problem for applied dynamical systems---
there was no  rigorous mathematical proof of his findings
until the 1999 paper of {\sc W. Tucker} \cite {Tucker99}.

The term ``chaos'' is used in many ways in  mathematical
literature.  In a widely accepted definition,   {\sc
  R. Devaney} \cite{Devaney86} defined a map  $f\co X\to X$
to be chaotic if it has these three  properties:
  \begin{quote} {\sc Dense periodic points:}  {\em Every nonempty open
   subset of $Y$ meets  a  periodic orbit of $f$.}
   
    {\sc Topological transitivity:}  {\em The orbit of some point of
      $Y$ is dense in $Y$.}

  {\sc Sensitivity to initial conditions:} {\em There exists $\del > 0$
     such that if $x, y \in Y$ are distinct, then  $\,d (f^kx - f^ky) \ge \del$
     \,for some $k>0$.}

  \end{quote}
In order to avoid trivial situations we add a fourth
property:\footnote
{Suggested by Devaney \cite{Devaney19}.}

\begin{quote} {\sc Nondiscreteness:} {\em Y is not a finite set},
\end{quote}
  When these hold I call  $Y$ and $f|Y$  {\em chaotic}.

{\sc P. Touhey} \cite{Touhey98}  proved  that  a remarkably
simple  condition is equivalent to Devaney's definition:
 \begin{quote} {\sc Sharing of periodic orbits:} {\em Every pair of nonempty
     open sets meet a common periodic orbit.}
 \end{quote}
 See also \cite{Touhey97, Touhey00}.
 
Devaney \cite[p. 324]{DHS04} validated Lorenz's conclusions by
constructing a Poincar\'e (or ``first return'' map)
 $f\co C\to   C$ for Lorenz's differential equations, having the
following properties:

 \begin{itemize}

\item   $C\subset \R 3$ is an affine open 2-cell.
  
\item  $f$ has a chaotic global attractor.
 \end{itemize} 
 
 It is easy to see that every periodic orbit fulfills the definition
 of ``chaotic,'' but  such  orbits are not of much dynamical
 interest.

 Our main result, Theorem 1, shows that
 monotone maps in strongly ordered spaces have no other chaotic
 attracting sets.
 
\subsection{Monotone dynamics}
The state space of a monotone systems is a space $X$ endowed with a
(partial) order, denoted by $\,\preceq$.  
The set $\{(x,y) \in X\times X\co x\preceq y\}$ is assumed to be closed. 

We write $x\prec y$ if $x\preceq y$ and $x\ne y$.  If $A$ and $B$ are
sets,
\[
A\prec B \dimply a\prec b, \quad (a\in A,\, b\in B).
\]
When $A$ is a singleton $\{a\}$ we also write $a \prec B$.

In the main result the ordered space $X$  is {\em strongly
  ordered}:  If $W\subset X$ is a neighborhood of $x$,  there are
nonempty open sets 
$U, V \subset W$ such that $U\prec  x \prec V$.

\smallskip
\noindent {\em Examples.}  Euclidean space $\R n$ is strongly ordered
by the classical {\em vector order}:
\[x\preceq y \dimply x_j\le y_j, \ (j=1,\dots,n).
\]
Many Banach spaces of real-valued functions   are strongly
ordered by the functional order, $f_1 \preceq f_2$ iff $f_1 x\le  f_2 x$ for all 
$x$ in the domain.  These spaces include the  spaces of $C^k$
functions on compact manifolds.

For $p >0$, $\mcal L_p$ spaces of
real-valued functions have the order  $f_1\preceq f_2$ iff $f_1 \le f_2$
almost everywhere.  These spaces are rarely  strongly ordered. 

A map $T$ between ordered spaces  is {\em monotone} if
\[x\preceq y  \implies Tx\preceq Ty.\]

In dynamical  models from   several  scientific fields,
including biology, chemistry, physics,
economics,  the set of states is an ordered space $S$, with the order
reflecting the relative ``size'' of states--- density, population, etc.
The evolution of states over time is modeled by a dynamical
system on $S$--- a family $F=\{F_t\}$  of maps
between subsets of $S$, closed under composition.  The time
variable $t$ varies over  either real numbers  or
integers.

Monotonicity means that the different species cooperate:  an increase in
the growth rate of one tends to increase the sizes of the others.
In many real-world settings this is
plausible. For example: sheep and grass cooperate, in that grass feeds
the sheep and sheep fertilize  grass.    

In many cases the maps $F_t$ are monotone for $t\ge 0$,  and $\Phi$ 
is called a monotone system. A typical example is a system of
differential equations in the positive orthant $\Rnp:= [0,\infty)^n$:
\[\ode{x_i}{t} = x_iG_i (x_1,\dots,x_n), x_i \ge 0, \quad
i=1,\dots,n,\]
modeling an ecology of $n$ interacting species. Here $x_i$ and $G_i$
are proxies for the size and the {\em per capita} growth rate of
species $i$. The state space is $\Rnp$ with the vector order:
\[x \preceq y \dimply x_i\le y_i \quad \text{\em for all $i$}. \]
Monotonicity is readily established when the partial derivatives of
the $G_i$ are continuous, and
\[ i \ne j \implies  \pde{G_i}{x_j} \ge 0.\]

If each species reproduces only once a year, the ecology is modeled by
a map $T\co \Rnp\to \Rnp$, and the dynamics is monotone provided the partial
derivatives of $T$ are continuous and nonnegative.  

Monotone dynamical systems   often 
permit  reliable predictions of behavior. In many case
it can be proved that  typical trajectories have predictable fates, such
tending  toward fixed points or periodic orbits.  See for example references
\cite {BasLemmens10,  DELEENHEER17,  DIRR, ENCISOSONTAG06,
Hirsch82,
Hirsch83,
  Hirsch88,
  Hirsch17,
  HS05,
 LajmanovichYorke,
  Landsberg96,
  Potsche15,
 Selgrade80,
  Smale76,
Smith95,
Zhao94}.
The recent survey by {\sc H. Smith} \cite {Smith17}  has an extensive
bibliography.

Monotonicity and chaos play quite different roles in applied dynamics.
Monotonicity is sometimes deliberately postulated, and can ordinarily be deduced
from the form of defining equations without extensive computations.
Monotonicity is useful because it usually leads to predictable
long-term behavior.

But chaos is undesirable: it makes accurate long-term
prediction is impossible, and is quite difficult to either prove or
disprove.  But as Lorenz discovered, simple models of realistic
systems can be chaotic.
\subsection{Results}
Our main result is  very simply stated:
\begin{quote}{\em A  monotone map  in a strongly  ordered
    space cannot  have a chaotic attractor.}
\end{quote}
In fact  slightly more is true:
\begin{theorem*} 
   Assume $X$ is strongly  ordered,   $\,T\co X\to X$  is monotone, 
 $\,A\subset X$ is attracting for $T$, and  one of the following
   conditions is satisfied:

  \begin{description}

  \item[(a)]  Periodic points are dense in $A$.

  \item[(b)] Some orbit  is dense in $A$.
  \end{description}
Then  $A$ is
a periodic orbit, and therefore  not chaotic. 
\end{theorem*}
\begin{proof}
  We rely on a {\em deus ex machina},
\cite[Theorem 4.1]{Hirsch85}:
 \begin{quote} { \em  Some  periodic orbit $O\subset A$  is
     attracting for $T\co X\to X$.}
 \end{quote}
We show that $A=O$.  Assume {\em per contra} $A\ne O$. 
Let  $W\subset X$ be the attractor basin of $O$.  Then
\begin{equation}\label{eq:2}
 x\in W \implies \lim_{k\to\infty}\dist \,(T^kx, O) =0.
\end{equation}
But (\ref{eq:2}) cannot be true if $A\ne O$: There exists $x\in W\cap
A\,\verb=\=O$ whose orbit is periodic when (a) holds, or dense in the
nonempty open set $A\sm O$ when (b) holds.
\end{proof}



\begin{thebibliography}{99}
  
\bibitem {BasLemmens10} M. Akian, S. Gaubert \& Bas Lemmens, {\em
   Stability and convergence in discrete convex monotone dynamical
   systems,} arXiv: 1003.5346v1 (2010).
 
 
 \bibitem {DELEENHEER17} P. De Leenheer, {\em The puzzle of partial
 migration}, Journal of Theroretical Biology {\bf 412} (2017),
   172--185.
  
 \bibitem{Devaney86}  R. Devaney, ``An introduction to chaotic
   dynamical systems,''  Benjamin/Cummings, Menlo Park, CA, (1986).

\bibitem{DHS04} R. Devaney, M. Hirsch \& S. Smale, `` Differential
  Equations, Dynamical Systems \& an Introduction to Chaos,'' Elsevier
  Academic Press 2004.

 \bibitem{Devaney19}  R. Devaney,  personal commuication (2019).
  
\bibitem{DIRR} G. Dirr et al., {\em Separable Lyapunov functions
for monotone systems: constructions and limitations,} Discrete and
Continuous Dynamical Systems,  Series B {\bf 20} (2015), 2497-–2526.
   
\bibitem{ENCISOSONTAG06} G. Enciso \& E. Sontag, {\em Global
  attractivity, i/o monotone small-gain theorems, and biological delay
  systems}, Discrete and Continuous Dynamical Systems {\bf 14} (2006),
  249--578

\bibitem{Hirsch85} M. Hirsch, {\em Attractors for discrete--time
  monotone dynamical systems in strongly ordered spaces.}  Geometry
  and Topology: Lecture Notes in Mathematics 1167,
  141--153. J. Alexander, J.Harer, editors.  Springer-Verlag, New
  York, 1985.

\bibitem{Hirsch82} M. Hirsch, Systems of differential
  equations which are competitive or cooperative. I: limit sets. SIAM
  J. Appl. Math. 13 (1982), 167--179.

\bibitem{Hirsch85a} M. Hirsch, Systems of differential equations
  which are competitive or cooperative, II: convergence almost
  everywhere, {\em SIAM J. Math. Anal.}, {\bf 16}, (1985) 423--439.
  
\bibitem{Hirsch88} M. Hirsch, {\em Stability and convergence in strongly
   monotone dynamical systems},  J. reine und angewandte
   Mathematik {\bf 383}  (1988),   1--53.

\bibitem{Hi88a} M. Hirsch, Systems of differential equations
  which are competitive or cooperative, III: competing
  species, {\em Nonlinearity {\bf 1}}, (1988) 51--71.
   
  \bibitem{Hirsch83} M. Hirsch, Differential equations and
    convergence almost everywhwere in strongly monotone
   semiflows, {\em Contemp. Math.}, {\bf 17}, (1983) 267--285.
 
 \bibitem{Hirsch84} M. Hirsch, The dynamical systems approach to
differential equations, Bull. Amer. Math. Soc. {\bf 11} (1984),
   1--64.

\bibitem{Hirsch17} M. Hirsch, {\em Monotone dynamical systems with
   polyhedral order cones and dense periodic points,} AIMS Mathematics
   {\bf 2} (2017), 24--27.  Also in: http://arxiv.org/abs/1611.09251
  
\bibitem{HS05} M. Hirsch \& H. Smith, {\em Monotone
  Dynamical Systems}, ``Handbook of Differential Equations,'' volume
  2, chapter 4.  A. Ca\~{n}ada, P.\ Dra\'{b}ek \& A.\ Fonda, editors.
  Elsevier North Holland 2005.

\bibitem{HSW94}  S.-B. Hsu,  H. Smith  \& P.  Waltman, {\em Dynamics  of
competition  in the unstirred chemostat}, Canadian
Appl. Math. Quart. {\bf 2} (1994), 461--483.  
  
\bibitem{Kolesov09} A. Kolesov \& N. Rozov, {\em On the definition of
  chaos,} Uspekhi Mat. Nauk {\bf 64} (2009), 125--172; translation in
Russian Math. Surveys {\bf 64} (2009), no. 4, 701–744 .

\bibitem{Kulczyncki08} M. Kulczyncki, {\em Noncontinuous maps and
  Devaney's chaos},   Regul. Chaotic Dyn. 13 (2008), no. 2, 81--84.
  
\bibitem{LajmanovichYorke} A. Lajmanovich \&  J. Yorke, {\em A
deterministic model for gonorrhea in a nonhomogeneous
population}, Mathematical  Biosciences {\bf 28}  (1976), 221--236.
  
\bibitem{Landsberg96} A. Landsberg \& E. Friedman, {\em Dynamical
    Effects of Partial Orderings in Physical Systems,}  Physical
    Review E {\bf 54} (1996), 3135--3141.
    


\bibitem{Lorenz63} E. Lorenz, {\em Deterministic Non-periodic Flow},
    J. Atmos. Sci. {\bf 20} (1963), 130--141.
       
\bibitem{MS90} J. Mallet-Paret \& H. Smith,
{\em The  Poincar\'e-Bendixson theorem for monotone cyclic
feedback systems}, J. Dynamics and  Diff. Equations
{\bf 2} (1990), 367--421. 
  
\bibitem{Potsche15} C. Potzsche, {\em Order-preserving nonautonomous discrete dynamics:
  Attractors and entire solutions}  Positivity {\bf 19} (2015), 547--576.

\bibitem{Selgrade80} J. Selgrade, {\em Asymptotic behavior of solutions to
single loop positive feedback systems}, J. Diff. Eqns. {\bf 38}
(1980), 80--103. 
  
\bibitem{Smale76} S. Smale {\em On the differential equations of species in
  competition,}   J. Math. Biol., {\bf 3} (1976),  5--7.


\bibitem{Smil84} J. Smillie, {\em Competitive and cooperative
tridiagonal systems of differential equations}, SIAM
J. Math. Anal. {\bf 15} (1984), 530-534. 

\bibitem{Sm86a} H. Smith, {\em Periodic competitive differential
equations and the discrete dynamics of competitive maps},
J. Diff. Eqns. {\bf 64} (1986), 165-194.

\bibitem{Smith95} H. Smith, ``Monotone Dynamical Systems: an
  introduction to the theory of competitive and cooperative systems,''
  Amer. Math. Soc. Surveys and Monographs {\bf 41}, 1995.

\bibitem{Smith17} H. Smith,  {\em Monotone dynamical systems: reflections on new
advances \& applications,} Discrete and Continuous Dynamical Systems
  Series B {\bf 37} (2017), no. 1, 485--504.

\bibitem{Touhey97} P. Touhey, {\em Yet another definition of chaos,} 
Amer. Math. Monthly {\bf 104} (1997), no. 5, 411–414.
     
\bibitem{Touhey98} P. Touhey, {Chaos: the evolution of a definition,}
Irish Math. Soc. Bull. No. 40 (1998), 60–70.

\bibitem{Touhey00} P. Touhey, {\em Persistent properties of chaos,} 
  J. Differ. Equations Appl. {\bf 6} (2000), no. 3,
  249–256.

   
\bibitem{Tucker99} W. Tucker,   {\em The Lorenz attractor exists},
  C. R. Acad. Sci. Paris, t. 328, S\'{e}rie I (1999), 1197--1202.
 
\bibitem {Zhao94} X.-Q. Zhao, {\em Global attractivity and stability
for discrete strongly monotone dynamical systems with applications to
biology}, Technical Report No. 94005.  Beijing: Inst. Appl. Math.,
Ac. Sinica (1994).  

\end{thebibliography}
\end{document}